\documentclass[secthm]{elsart}

\usepackage{latexsym}  

 \makeatletter     

 \def\@journal{Discrete Mathematics}

\newtheorem{olddefn}[thm]{Definition}
\renewenvironment{defn}{\begin{olddefn}\upshape}{\end{olddefn}}
\newtheorem{oldassum}[thm]{Assumption}
\renewenvironment{assum}{\begin{oldassum}\upshape}{\end{oldassum}}

\renewcommand\Elproofname{Proof.}  

\renewenvironment{pf}%
  {\par\noindent
   {\bfseries\Elproofname}\enspace\ignorespaces}%
  {\par}
\expandafter\let\csname endpf*\endcsname=\endpf 

  \def\@opargbegintheorem#1#2#3{\trivlist
     \item[\hskip\labelsep{\bfseries #1\ #2}\ {(#3).}]\itshape}
  \def\@begintheorem#1#2{\trivlist
    \item[\hskip \labelsep{\bfseries #1\ #2.}]\itshape}

\def\[{\relax\ifmmode\@badmath
  \else
  \@beginparpenalty\predisplaypenalty
  \@endparpenalty\postdisplaypenalty
  \begin{trivlist}\@topsep \eqntopsep       
   \@topsepadd -6pt                   
  \item[]\leavevmode
   \hbox to\linewidth\bgroup\hfil $ \displaystyle
  \hskip\mathindent\bgroup\fi}
\def\]{\relax\ifmmode \egroup $\hfil \egroup
  \end{trivlist}
  \else \@badmath \fi}
\let\oldeqnarray=\eqnarray
\def\eqnarray{\vspace{-12pt}\oldeqnarray}

 \@latex@warning@no@line
 {Reminder -- we need to delete this section\MessageBreak
 before sending the TeX file to Elsevier} 

 \makeatother     

\newcommand\Cay{\operatorname{Cay}}

\begin{document}

\begin{frontmatter}

\title{Automorphism Groups with Cyclic Commutator Subgroup
 and Hamilton Cycles}

\author{Edward Dobson}
 \address{Department of Mathematics,
 Louisiana State University,
 Baton Rouge, LA 70803, USA}

\author{Heather Gavlas}
 \address{Department of Mathematics and Statistics,
 Grand Valley State University,
 Allendale, MI 49401, USA}

\author{Joy Morris}
 \address{Department of Mathematics and Statistics,
 Simon Fraser University,
 Burnaby, BC V5A~1S6, Canada}

\author{Dave Witte}
 \address{Department of Mathematics,
 Oklahoma State University,
 Stillwater, OK 74078, USA}

\begin{abstract}
 It has been shown that there is a Hamilton cycle in every connected Cayley
graph on each group $G$ whose commutator subgroup is cyclic of prime-power
order. This paper considers connected, vertex-transitive graphs~$X$ of order
at least~$3$ where the automorphism group of~$X$ contains a transitive
subgroup~$G$ whose commutator subgroup is cyclic of prime-power order. We
show that of these graphs, only the Petersen graph is not hamiltonian.
 \end{abstract}


 \end{frontmatter}

\section{Introduction}
 Considerable attention has been devoted to the problem of determining
whether or not a connected, vertex-transitive graph~$X$ has a Hamilton cycle
\cite{AlspachSurvey}, \cite{CurranGallian}, \cite{WitteGallian}. A graph $X$
is \emph{vertex-transitive} if some group~$G$ of automorphisms of~$X$ acts
transitively on~$V(X)$. If $G$ is abelian, then it is easy to see that $X$
has a Hamilton cycle. Thus it is natural to try to prove the same conclusion
when $G$ is ``almost abelian." Recalling that the \emph{commutator subgroup}
of~$G$ is the subgroup $G' = \langle x^{-1}y^{-1}xy:x,y\in G\rangle$, and that
$G$ is abelian if and only if the commutator subgroup of~$G$ is trivial, it
is natural to consider the case where the commutator subgroup of~$G$ is
``small" in some sense. In this vein, K.~Keating and D.~Witte
\cite{KeatingWitte} used a method of D.~Maru\v si\v c \cite{Marusic} to show
that there is a Hamilton cycle in every Cayley graph on each group whose
commutator subgroup is cyclic of prime-power order. This paper utilizes
techniques of B.~Alspach, E.~Durnberger, and T.~Parsons
\cite{AlspachMetaHam}, \cite{AlspachMetaHamPrime},
\cite{AlspachMetaHamPower} to prove the following result.

\begin{thm} \label{mainthm}
 Let $X$ be a connected vertex-transitive graph of order at least $3$. If
there is a transitive group~$G$ of automorphisms of~$X$ such that the
commutator subgroup of~$G$ is cyclic of prime-power order, then $X$ is the
Petersen graph or $X$ is hamiltonian.
 \end{thm}

Because $K_2$ and the Petersen graph have Hamilton paths, the following
corollary is immediate.

\begin{cor}
 Let $X$ be a connected vertex-transitive graph. If there is a transitive
group~$G$ of automorphisms of~$X$ such that the commutator subgroup of~$G$
is cyclic of prime-power order, then $X$ has a Hamilton path.
 \end{cor}

\section{Assumptions and Definitions} \label{assumps-defs}

\begin{assum} \label{BasicAssumps} Throughout this note, $X$ is a connected
vertex-transitive graph, $G$ is a group of
automorphisms of~$X$ that acts transitively on the vertex set $V(X)$, and
$G'$ is the commutator subgroup of $G$.
 \end{assum}

Although the following definitions and results may be stated in more general
group-theoretic terms (see \cite{Scott} or~\cite{Biggs}), we state them here
in the context of this problem.

\begin{defn}
 The \emph{stabilizer} $G_x$ of a vertex $x\in V(X)$ is $\{\,g \in G : g(x) =
x \,\}$ and is a subgroup of~$G$.
 \end{defn}

\begin{lem}[{\cite[10.1.2, p.~256]{Scott}}] \label{StabConj}
 Let $x \in V(X)$ and $g \in G$. Then $G_{gx} = g(G_x)g^{-1}$.
 \end{lem}

\begin{cor} \label{StabNorm}
 If $H$ is a normal subgroup of~$G$, then the following are equivalent:
 \begin{enumerate}
 \item \label{HGx-normal-some}
 $H G_x$ is a normal subgroup of~$G$ for some $x \in V(X)$;
 \item \label{HGx-normal-every}
 $H G_x$ is a normal subgroup of~$G$ for every $x \in V(X)$;
 \item \label{HGx=Hgy}
 $H G_x = H G_y$ for all $x,y \in V(X)$.
 \end{enumerate}
 \end{cor}

\begin{pf} Let $g \in G$ and $x \in V(X)$. From the lemma, we know $G_{gx} =
g(G_x) g^{-1}$, and since $H$ is normal, we have $H = gHg^{-1}$. So
 \begin{equation} \label{HGx}
 H G_{gx} = (gHg^{-1})(g(G_x) g^{-1}) = g (HG_x)g^{-1} .
 \end{equation}

 (\ref{HGx-normal-some})~$\Rightarrow$~(\ref{HGx=Hgy}).
 Since $G$ is transitive on $V(X)$, there exists $g \in G$ with $gx = y$.
Then, since $HG_x$ is normal, (\ref{HGx})~implies $HG_y = HG_x$, as desired.

 (\ref{HGx=Hgy})~$\Rightarrow$~(\ref{HGx-normal-every}).
 Let $g \in G$. From~(\ref{HGx=Hgy}), we have $H G_{gx} = H G_x$. Therefore,
(\ref{HGx})~implies $g(HG_x)g^{-1} = HG_x$, as desired. \qed
 \end{pf}

\begin{cor} \label{NoNormalInStab}
 For every $x \in V(X)$, the stabilizer $G_x$ does not contain a nontrivial,
normal subgroup of~$G$.
 \end{cor}

\begin{pf} Let $H$ be a normal subgroup of~$G$ that is contained in~$G_x$.
Lemma~\ref{StabConj} implies $H \subset G_{gx}$, for all $g \in G$. Since
$G$ is acts transitively on $V(X)$, it follows that $H \subset G_y$, for all
$y \in V(X)$. Therefore, the identity automorphism of~$X$ is the only
element of~$H$. \qed
 \end{pf}

\begin{defn} Let $H$ be a subgroup of~$G$, and let $x \in V(X)$. The 
\emph{$H$-orbit} of~$x$ is $\{\, hx : h \in H \,\}$. The $H$-orbits form a
partition of~$V(X)$, and if $H$ is normal in~$G$, then the subgraphs of~$X$
induced by distinct $H$-orbits are isomorphic, as $g(Hx) = H(gx)$ in this
case.
 \end{defn}

\begin{defn} Let $H$ be a subgroup of~$G$. The \emph{quotient graph} $X/H$
is that graph whose vertices are the $H$-orbits, and two such vertices $Hx$
and $Hy$ are adjacent in $X/H$ if and only if there is an edge in $X$
joining a vertex of $Hx$ to a vertex of~$Hy$. If $H$ is normal in~$G$, then
the action of~$G$ on~$V(X)$ factors through to a transitive action of $G/H$
on $V(X/H)$ by automorphisms of~$X/H$ and thus $X/H$ is vertex-transitive.
 \end{defn}

\begin{lem} \label{LiftPath}
 If $H$ is a normal subgroup of~$G$, then every path in $X/H$ lifts to a
path in~$X$.
 \end{lem}

\begin{pf} It suffices to show that if $Hx$ is adjacent to~$Hy$ in $X/H$,
then $x$~is adjacent to some vertex in~$Hy$. By definition of $X/H$, we know
that some $\tilde x \in Hx$ is adjacent to some $\tilde y \in Hy$. Next
there exists $h \in H$ with $x = h \tilde x$, so that $x$ is adjacent to $h
\tilde y \in Hy$. \qed
 \end{pf}

\begin{defn} Let $S$ be a subset of~$G$, and assume $S$ is \emph{symmetric}
(that is, $s^{-1} \in S$ for all $s \in S$). The \emph{Cayley graph}
$\Cay(G;S)$ is that graph whose vertices are the elements of~$G$, and for
vertices $g$ and $h$, there is an edge from~$g$ to~$h$ if and only if $gs
=h$ for some $s \in S$. Since $G$ acts transitively on the vertices
of~$\Cay(G;S)$ by left multiplication, $\Cay(G;S)$ is vertex-transitive. A
Cayley graph is connected if and only if $S$ generates~$G$.
 \end{defn}

Recall that $G'$ is a normal subgroup of~$G$ and that the quotient group
$G/G'$ is abelian \cite[Thms.~3.4.11 and 3.4.10, p.~59]{Scott}. Since $G/G'$
is abelian and transitive on $V(X/G')$, it follows from the next result that
$X/G'$ is a Cayley graph on the abelian group $G/(G_x G')$, for any $x \in
V(X)$.

\begin{lem}[Sabidussi \cite{Sabidussi}] \label{SabThm}
 If $G_x$ is trivial for some $x \in V(X)$, then $X$ is (isomorphic to) a
Cayley graph on~$G$.
 \end{lem}

\section{Preliminaries on the Frattini subgroup}

As in Section~\ref{assumps-defs}, we assume that
Assumption~\ref{BasicAssumps} holds. 

\begin{assum} \label{assume-cyclic}
 We assume $G'$ is cyclic of order~$p^k$, where $p$~is a prime, and that $X$
has at least three vertices.
 \end{assum}

\begin{assum} \label{Gminimal}
 We also assume $X$ is $G$-minimal. That is, if $Y$ is a connected, spanning
subgraph of~$X$, such that, for all $g \in G$, we have $gY = Y$, then it
must be the case that $Y = X$. (In the case of Cayley graphs, $\Cay(G;S)$ is
$G$-minimal if and only if no proper symmetric subset of~$S$ generates~$G$.)
Since a Hamilton cycle in any such subgraph~$Y$ would also be a Hamilton
cycle in~$X$, we may assume this without loss of generality.
 \end{assum}

The main result of this section is Lemma~\ref{edgeHamiltonian}. A central
idea to the proof is that of the Frattini subgroup, defined in
\cite[\S7.3]{Scott}.

\begin{defn} An element~$g$ of~$G$ is a \emph{nongenerator} if, for every
subset~$S$ of~$G$ such that $\langle S,g\rangle = G$, we have that $\langle S
\rangle = G$.
 The \emph{Frattini subgroup} of~$G$, denoted~$\Phi(G)$, is the set of all
nongenerators of~$G$ and is a subgroup of~$G$.
 \end{defn}

\begin{lem} \label{HpFrattini}
 If $H$ is any subgroup of~$G'$, then $H$ is normal in~$G$ and $H^p \subset
\Phi(G)$, where
 $H^p = \langle h^p : h \in H \rangle$.
 \end{lem}

\begin{pf} Since $G'$ is a cyclic normal subgroup of~$G$, we know that every
subgroup of~$G'$ is a normal subgroup of~$G$ \cite[Thm.~1.3.1(i), p.~9, and
Thm.~2.1.2(ii), p.~16]{Gorenstein}. Therefore $H$ is normal in~$G$ and hence
$\Phi(H) \subset \Phi(G)$ \cite[7.3.17, p.~162]{Scott}. Since $H$ is a
cyclic $p$-group, it is not difficult to see that $\Phi(H) = H^p$
\cite[7.3.7, p.~160]{Scott}. \qed
 \end{pf}

\begin{lem} \label{Frattini}
 If $H$ is a normal subgroup of~$G$ and $H \subset \Phi(G)$, then $X/H$ is
$G$-minimal.
 \end{lem}

\begin{pf} Let $Y$ be a connected, spanning subgraph of~$X/H$ such that for
all $g \in G$, we have that $gY = Y$. Let $x \in V(X)$, and let
 \begin{eqnarray*}
 S &=& \{\, s \in G : \hbox{$sx$ is adjacent to $x$ in $X$} \,\},
 \hbox{ and} \\
 T &=& \{\, t \in G : \hbox{$Htx$ is adjacent to $Hx$ in $Y$} \,\} .
 \end{eqnarray*}
 It is straightforward to verify that $G_x S G_x = S$ and $HG_x T G_x = T$.
Furthermore, since $Y$ is connected, we see that $T$ generates~$G$.

Since $Y$ is a subgraph of~$X/H$, it must be the case that $T \subset HS$.
Hence, since $HT = T$, we have that $T = H(S \cap T)$. Next since $T$
generates~$G$ and $H \subset \Phi(G)$, we conclude that $S \cap T$
generates~$G$. Therefore, letting $Z$ be the spanning subgraph of~$X$ whose
edge set is
 \[ E(Z) = \Bigl\{\, \{gtx, gx\} : g \in G, ~ t \in S \cap T \,\Bigr\}, \]
 we see that $Z$~is connected. So $Z$ is a connected, spanning subgraph of
$X$ such that $gZ = Z$ for all $g \in G$. Since $X$ is $G$-minimal, it
follows that $Z = X$ and hence $S \cap T = S$. Therefore $HS = H(S \cap T) =
T$, so $X/H = Y$. \qed
 \end{pf}

Because a $G$-minimal graph has no loops, we have the following corollary.

\begin{cor} \label{noloops}
 If $H$ is a normal subgroup of~$G$ and $H \subset \Phi(G)$, then the
subgraph of~$X$ induced by each $H$-orbit has no edges.
 \end{cor}

We now recall (in a weak form) the fundamental work of C.~C.~Chen and
N.~F.~Quimpo \cite{ChenQuimpo}.

\begin{thm}[Chen-Quimpo \cite{ChenQuimpo}] \label{ChenQuimpoThm}
 Let $Y$ be a connected Cayley graph on an abelian group of order at least
three. Then each edge of~$Y$ (except any loop) is contained in some Hamilton
cycle of~$Y$.
 \end{thm}

The following helpful result is the main conclusion obtained from our
discussion of $G$-minimality and Frattini subgroups. (It also relies on the
Chen-Quimpo Theorem.)

\begin{lem} \label{edgeHamiltonian}
 If $H$ is a subgroup of $G'$ such that $X/H$ has a Hamilton cycle, then
each edge of~$X/H$ (except any loop) is contained in some Hamilton cycle
of~$X/H$.
 \end{lem}

\begin{pf} If $H = G'$, then we have already seen that $X/G'$ is a Cayley 
graph on the abelian group $G/(G_x G')$ and hence desired conclusion follows
from the Chen-Quimpo Theorem (\ref{ChenQuimpoThm}).

We may now assume $H \not=G'$, which implies $H \subset (G')^p$. So $H
\subset \Phi(G)$ by Lemma~\ref{HpFrattini}; therefore $X/H$ is $G$-minimal
by Lemma~\ref{Frattini}. Let $C$ be a Hamilton cycle in~$X/H$, and let $Y =
\cup_{g \in G}\ gC$. Since $X/H$ is $G$-minimal, we must have $Y = X/H$, and
thus every edge of $X/H$ is contained in some Hamilton cycle~$gC$. \qed
 \end{pf}

\section{Proof of Theorem~\ref{mainthm}}

As before, we assume that Assumptions \ref{BasicAssumps},
\ref{assume-cyclic}, and~\ref{Gminimal} still hold. The main conclusions of
this section are  two propositions which together constitute a proof of
Theorem~\ref{mainthm}.

Let us begin by disposing of a trivial case, namely the case when $X/G'$ has
only one vertex. Then $G'$ is transitive on~$V(X)$. Furthermore, we see from
Corollary~\ref{NoNormalInStab} and Lemma~\ref{HpFrattini} that $(G')_x =
\{e\}$ for each vertex~$x$ of $X$. Thus it follows by Lemma~\ref{SabThm} that
$X$ is a Cayley graph on the abelian group $G'$. Then
Theorem~\ref{ChenQuimpoThm} implies that $X$ has a Hamilton cycle if $X$ has
order at least~3.

\begin{lem} \label{FactorGroupLemma}
 Suppose $H$ is a subgroup of~$G'$ and that
 \[ H^p x_1, H^p x_2, \ldots, H^p x_n, H^p x_{n+1} \]
 is a path in~$X/H^p$ with $H^px_1 \not= H^p x_{n+1}$. If $H x_1, H x_2,
\ldots, H x_n, H x_{n+1}$ is a Hamilton cycle in $X/H$ (or if we have $n =
2$, $X/H \cong K_2$, and $H x_1 = H x_3 \not= H x_2$), then $X$ has a
Hamilton cycle.
 \end{lem}

\begin{pf} By Lemma~\ref{LiftPath}, we can lift the path $H^p x_1, H^p x_2, 
\ldots, H^p x_n, H^p x_{n+1}$ in $X/H^p$ to a path $x_1,x_2,\ldots,x_{n+1}$ 
in~$X$. Since $H x_1 = H x_{n+1}$, there exists $\gamma \in H$ such that 
$\gamma(x_1) = x_{n+1}$. Now, since $x_{n+1} \not\in H^p x_1$, it follows
that $\gamma \not\in H^p$, which implies that $\gamma$ generates~$H$. Let
$P$ be the path $x_1,x_2,\ldots,x_n$. Then the trail 
$P,\gamma(P),\ldots\gamma^{|H|-1}(P),x_1$ is a Hamilton cycle in $X$. \qed
 \end{pf}

The analysis now breaks into two cases, depending on whether the subgraphs
induced by each $G'$-orbit are empty. Since $G'$ is a normal subgroup, all
of these subgraphs are isomorphic, and hence either all are empty, or none
are.

\begin{prop} \label{noedges}
 If the subgraph induced by each $G'$-orbit is empty, then $X$ has a
Hamilton cycle.
 \end{prop}

\begin{pf*}{Proof \rm (cf.~\cite{AlspachMetaHam}, \cite{AlspachMetaHamPrime},
\cite{AlspachMetaHamPower}).}
 Let $x_1 \in V(X)$. Since $G/G'$ is abelian, it follows that $G' G_{x_1}$
is a normal subgroup of~$G$. Hence, there is a subgroup~$H$ of~$G'$, such
that $H G_{x_1}$ is normal in~$G$, but $K G_{x_1}$ is \emph{not} normal
in~$G$, for every proper subgroup~$K$ of~$H$. (It may be the case that $H =
G'$ or $H = \{e\}$.) Since $X/H$ is a connected Cayley graph on the group
$G/(H G_x)$ (see Lemma~\ref{SabThm}) and the commutator subgroup of $G/HG_x$
is cyclic, it follows that $X/H$ has a Hamilton cycle or $X/H \cong K_2$
\cite{KeatingWitte}.

We may assume that $H \not= \{e\}$, for otherwise $X = X/H$ has a Hamilton
cycle, and we are done. Then $H^p \not= H$, and the choice of~$H$ implies
that $H^p G_{x_1}$ is not normal in $G$. Therefore, since $X$ is connected
and vertex-transitive, it follows from Corollary~\ref{StabNorm} that $x_1$ is
adjacent to some vertex~$u$ such that $H^p G_{x_1} \not= H^p G_u$. This
implies that there exists $\gamma\in G_{x_1}$ such that $\gamma(u)\not \in
H^p u$. However, since $H G_{x_1} = H G_u$ (see Corollary~\ref{StabNorm}),
we have that $\gamma(u) \in G_{x_1} u \subset H G_u u = Hu$.

Since the subgraph induced by~$H x_1$ is contained in the subgraph induced
by $G'x_1$, which has no edges, and $x_1$ is adjacent to~$u$, it follows
that $u\not\in H x_1$, and thus $\{Hx_1,Hu\}$ is an edge in~$X/H$.
Therefore, there exists a Hamilton path from $H x_1$ to $H u$ in $X/H$ (see
Lemma~\ref{edgeHamiltonian}). This path lifts to a path $x_1,x_2,x_3,\ldots,
x_n$ in~$X$, where $x_n\in H u$ (see Lemma~\ref{LiftPath}). Since not both of
 \[ H^p u,H^p x_1,H^p x_2,\ldots,H^p x_n
 \hbox{\qquad and\qquad}
 H^p \gamma(u),H^p x_1,H^p x_2,\ldots,H^p x_n \]
 can be a cycle, Lemma~\ref{FactorGroupLemma} implies there is a Hamilton
cycle in $X$ as desired. \qed
 \end{pf*}

We now consider the case where the $G'$-orbits do not induce empty graphs.
Let us begin with some preliminary observations.

\begin{lem} \label{connectedAndOdd}
 If each subgraph induced by each $G'$-orbit is nonempty, then these
subgraphs are connected and $p$~is odd.
 \end{lem}

\begin{pf} Suppose that the subgraph induced by~$G'x$ is not connected.
Since $G'$ is cyclic, this subgraph is circulant, and hence each connected
component must be induced by the orbit of some proper subgroup~$H$ of~$G'$.
But $H \subset (G')^p$, and $(G')^p \subset \Phi(G)$ (see
Lemma~\ref{HpFrattini}), and Corollary~\ref{noloops} asserts that the
subgraph induced by any $H$-orbit has no edges. This contradicts the fact
that the connected components of the subgraph induced by $G'x$ \emph{do}
have edges.

We now show that $p$ is odd. Suppose, to the contrary, that $p = 2$. Let
$\bar{G} = G/(G')^2$. The commutator subgroup of~$\bar{G}$ is $G'/(G')^2$,
which has order~2. Because a group of order~2 has no nontrivial
automorphisms, this implies that the commutator subgroup of~$\bar{G}$ is
contained in the center of~$\bar{G}$; therefore $\bar{G}$ is nilpotent (of
class~2) \cite[p.~21]{Gorenstein}. Since $(G')^2 \subset \Phi(G)$ (see
Lemma~\ref{HpFrattini}), it follows that $G/\Phi(G)$ is nilpotent. Hence $G$
itself is nilpotent \cite[7.4.10, p.~168]{Scott}, so $G' \subset \Phi(G)$
\cite[Thm.~7.3.4, p.~160]{Scott}. Therefore the subgraph induced by each
$G'$-orbit is empty (see Corollary~\ref{noloops}), contradicting our
hypothesis.\qed
 \end{pf}

We can now concisely state several important results of
B.~Alspach~\cite{AlspachMetaHamPower}, \cite{AlspachLifting}. They have been
rephrased in the context of our problem.

\begin{thm}[Alspach] \label{AlspachThm}
 Assume that the subgraph induced by each $G'$-orbit is nonempty. Then $X$
has a Hamilton cycle if any of the following are true:
 \begin{enumerate}
 \item \label{AlspachThmA}
 every vertex of the subgraph induced by a $G'$-orbit has degree at
least~$3$\/ {\rm \cite[Thm.~2.4]{AlspachLifting}}; or
 \item \label{AlspachThmB}
 $X/G'$ has only two vertices and $X$ is not the Petersen graph\/ {\rm
\cite[Thm.~2]{AlspachMetaHamPower}}; or
 \item \label{AlspachThmC}
 the number of vertices of $X/G'$ is odd\/ {\rm
\cite[Thm.~3.7(ii)]{AlspachLifting}}; or
 \item \label{AlspachThmD}
 there is a Hamilton cycle in $X/G'$ that can be lifted to a cycle in~$X$\/
{\rm \cite[Thm.~3.9]{AlspachLifting}}.
 \end{enumerate}
 \end{thm}

\begin{lem} \label{StabsAllSame}
 Let $x \in V(X)$. If $G_x = G_y$ for all $y \in G'x$, then $X$ has a
Hamilton cycle.
 \end{lem}

\begin{pf} This is essentially the same as the proof of
Proposition~\ref{noedges}; the assumption that the subgraph induced by~$G'x$
has no edges was used only to show that $u \not\in H x_1$, and this follows
from the assumption that $G_x = G_y$ for all $y \in G'x$ (and hence for all
$y \in Hx$). \qed
 \end{pf}

The following lemma shows that we may assume that all the vertices in each
$G'$-orbit have different stabilizers. The proof is mainly group-theoretic.
The key observation is that the automorphism group of a cycle is a dihedral
group. Therefore, if a group of automorphisms acts transitively on the
vertices of an odd cycle, then either all vertices have different
stabilizers or all vertices have the same stabilizer, depending on whether
 the group contains a reflection.

\begin{lem} \label{SomeStabsSame}
 Assume that the subgraph induced by each $G'$-orbit is nonempty, and that
there are two vertices $x$~and~$y$ belonging to the same $G'$-orbit such
that $G_x = G_y$. Then $X$ has a Hamilton cycle.
 \end{lem}

\begin{pf} Let $Y$ be the subgraph of~$X$ induced by $G' x$, and let $K =
\cap_{v \in G'x} G_v$. (Note that $K$ is a subgroup.) Since every subgroup
of~$G'$ is normal in~$G$ (see Lemma~\ref{HpFrattini}), it follows that $G' 
\cap G_x = \{e\}$ (see Corollary~\ref{NoNormalInStab}) and hence $G' \cap K
= \{e\}$. On the other hand, since $G'$ fixes $V(Y)$ setwise, we see from
Lemma~\ref{StabConj} that $G'$ normalizes~$K$. Therefore, $[G',K] \subset G'
\cap K$, so $G'$ must centralize~$K$.

By Theorem~\ref{AlspachThm}(\ref{AlspachThmA}), if every vertex of $Y$ has
degree at least 3, then $X$ has a Hamilton cycle. Thus we may assume that
$Y$ is 2-regular. Since $Y$ is connected and has an odd number of vertices
(see Lemma~\ref{connectedAndOdd}), it follows that $Y$ is a odd cycle.
Therefore, we see that $K$ is a subgroup of index at most two in~$G_v$, for
each $v \in V(Y)$. In fact, from Lemma~\ref{StabsAllSame}, we may assume
that the index is exactly two.

Let $A$ be a subgroup of $G_x$ of order two. Since $A$ is not normal in~$G$
(see Corollary~\ref{NoNormalInStab}), we know that $A$ does not
centralize~$G'$ (otherwise, it would be the only Sylow $2$-subgroup of the
normal subgroup $AG'$, and hence $A$ would be normal in~$G$). Since $G'$ is
a cyclic $p$-group and $p$~is odd, the automorphism group of~$G'$ is cyclic
\cite[5.7.12, p.~120]{Scott} and therefore has exactly one element of
order~2, namely, inversion. Therefore, the action of~$A$ by conjugation
inverts~$G'$. Since $G'$ has odd order, this means that $e$~is the only
element of~$G'$ that is centralized by~$A$.

On the other hand, $A$ must centralize~$K$ (since $A \subset G_x$, $G_x$
normalizes~$K$, and $K \cap G' = \{e\}$). Thus, we see that $K$ is the
centralizer of $AG'$ in $KG'$. Since $AG'$ and $KG'$ are normal, we have 
that $K$ is a normal subgroup of~$G$. Therefore, $K = \{e\}$ (see
Corollary~\ref{NoNormalInStab}), which implies $G_x = A$ has order~2. Hence,
since a group of order~2 has no nontrivial automorphisms, any element of~$G$
that normalizes $G_x$ must actually centralize it. In particular, then the
conclusion of the preceding paragraph implies that no nontrivial element
of~$G'$ normalizes $G_x$. This contradicts the fact that $G_x = G_y$ (see
Lemma~\ref{StabConj}). \qed
 \end{pf}

\begin{prop} \label{ifedges}
 If the subgraph induced by each $G'$-orbit has some edges, then $X$ has a
Hamilton cycle or $X$ is the Petersen graph.
 \end{prop}

\begin{pf*}{Proof \rm (cf.~pf.~of Prop.~\ref{noedges}).} Let $H$ be the
smallest subgroup of~$G'$ such that whenever $x$~and~$y$ are two adjacent
vertices of~$X$ not belonging to the same $G'$-orbit, we have $HG_x = HG_y$.
(It may be the case that $H = G'$.) Note that, from
Theorem~\ref{AlspachThm}(\ref{AlspachThmB}), we may assume $X/G'$ has more
than two vertices.

Assume for the moment that $H$ is nontrivial. Then $H^p$ is properly
contained in~$H$, so the minimality of~$H$ implies there are two adjacent
vertices $x_1$~and~$u$, such that $G' x_1 \not= G'u$, and $H^p G_{x_1} \not=
H^p G_u$. Thus, there exists $\gamma\in G_{x_1}$ such that $\gamma(u)\not
\in H^p u$. Since $X/G'$ has more than two vertices, we have that $X/H$ is
not the Petersen graph, and from Lemma~\ref{edgeHamiltonian} (and induction
on the number of vertices in~$X$), we know there is a Hamilton path
from~$Hx_1$ to~$Hu$ in $X/H$. This path lifts to a path $x_1,x_2,\ldots,x_n$
in~$X$, where $x_n \in Hu$ (see Lemma~\ref{LiftPath}). Since not both of
 \[ H^p u,H^p x_1,H^p x_2,\ldots,H^p x_n
 \hbox{\qquad and\qquad}
 H^p \gamma(u),H^p x_1,H^p x_2,\ldots,H^p x_n \]
 can be a cycle, Lemma~\ref{FactorGroupLemma} implies there is a Hamilton
cycle in $X$, as desired.

We may now assume $H = \{e\}$. Let $x_1, x_2, \ldots, x_{m+1}$ be a lift
in~$X$ of a Hamilton cycle in $X/G'$. Because $H = \{e\}$, we must have
$G_{x_{i}} = G_{x_{i+1}}$ for every~$i$, so $G_{x_1} = G_{x_{m+1}}$.
Therefore, if $x_{1} \not= x_{m+1}$, then Lemma~\ref{SomeStabsSame} implies
that $X$ has a Hamilton cycle. On the other hand, if $x_{1} = x_{m+1}$, then
Theorem~\ref{AlspachThm}(\ref{AlspachThmD}) yields the same conclusion. \qed
 \end{pf*}

\begin{ack} Much of this research was carried out at the Centre de Recherches
Math\'e\-matiques of the Universit\'e de Montr\'eal. The authors would like
to thank the organizers and participants of the Workshop on Graph Symmetry,
and the staff of the CRM, for the stimulating environment they provided. In
particular, they are grateful to Brian Alspach for his helpful suggestions
and encouragement, both during the course of this research and at other
times in their careers. Witte was partially supported by a grant from the
National Science Foundation.
 \end{ack}


\begin{thebibliography}{99}

\bibitem{AlspachSurvey}
 B. Alspach,
 The search for long paths and cycles in vertex-transitive graphs and
digraphs,
 in: K.L. McAvaney, ed., \emph{Combinatorial Mathematics VIII},
 Lecture Notes in Mathematics, Vol. 884
 (Springer-Verlag, Berlin, 1981) 14--22.

\bibitem{AlspachMetaHamPower}
 B. Alspach,
 Hamilton cycles in metacirculant graphs with prime power cardinal blocks,
 \emph{Ann. Discrete Math} {\bfseries 41} (1989) 7--16.

\bibitem{AlspachLifting}
 B. Alspach,
 Lifting Hamilton cycles of quotient graphs,
 \emph{Discrete Math.} {\bfseries 78} (1989) 25--36.

\bibitem{AlspachMetaHamPrime}
 B.~Alspach, E.~Durnberger, and T.~Parsons,
 Hamilton cycles in metacirculant graphs with prime cardinality blocks,
 \emph{Ann. Discrete Math} {\bfseries 27} (1985) 27--34.

\bibitem{AlspachMetaHam}
 B. Alspach and T.~Parsons,
 On hamiltonian cycles in metacirculant graphs,
 \emph{Ann. Discrete Math} {\bfseries 15} (1982) 1--7.

\bibitem{Biggs}
 N.~Biggs,
 \emph{Algebraic Graph Theory, 2nd ed.}
 (Cambridge Univ. Press, Cambridge, 1993).

\bibitem{ChenQuimpo}
 C.C. Chen and N.F. Quimpo,
 On strongly hamiltonian abelian group graphs,
 in: K.L. McAvaney, ed., \emph{Combinatorial Mathematics VIII},
 Lecture Notes in Mathematics, Vol. 884
 (Springer-Verlag, Berlin, 1981) 23--34.

\bibitem{CurranGallian}
 S.J.~Curran and J.A.~Gallian,
 Hamiltonian cycles and paths in Cayley graphs and digraphs --- a survey,
 \emph{Discrete Math.} {\bfseries 156} (1996) 1--18.

\bibitem{Gorenstein}
 D.~Gorenstein,
 \emph{Finite Groups}
 (Chelsea, New York, 1980).

\bibitem{KeatingWitte}
 K. Keating and D. Witte,
 On Hamilton cycles in Cayley graphs in groups with cyclic commutator
subgroup,
 \emph{Ann. Discrete Math} {\bfseries 27} (1985) 89--102.

\bibitem{Marusic}
 D. Maru\v si\v c,
 Hamiltonian circuits in Cayley graphs,
 \emph{Discrete Math} {\bfseries 46} (1983) 49--54.

\bibitem{Sabidussi}
 G. Sabidussi, Vertex-transitive graphs,
 \emph{Monatshefte fur Math.}, {\bfseries 68} (1964) 426--438.

\bibitem{Scott}
 W.~R.~Scott,
 \emph{Group Theory}
 (Dover, New York, 1987).

\bibitem{WitteGallian}
 D. Witte and J. A. Gallian,
 A survey: hamiltonian cycles in Cayley digraphs,
 \emph{Discrete Math.}, {\bfseries 51} (1984) 293--304.

\end{thebibliography}
\end{document}